\documentstyle[12pt]{article}
\begin{document}
\begin{center}
{\large {\bf TWO--PARAMETRIC DEFORMATION U$_{p,q}$[gl(2/1)]\\[2mm]
AND ITS INDUCED REPRESENTATIONS}}\\[5mm]
{\it In memory of Professor Asim Orhan Barut}\\[1mm]
\vskip 0.9truecm
\def\thefootnote{a)}
{\bf Nguyen Anh Ky}\hspace*{1.5mm}\footnote{~On leave of absence from the
Institute of Physics, National Centre for Natural Science and 
Technology, P.O. Box 429, Bo Ho, Hanoi 10000, {\bf Vietnam}}   
\def\thefootnote{ }
\\[1mm] 
\normalsize 
Theory Division, CERN, Geneve 23, CH--1211, {\bf Switzerland}
\vskip 1.7truecm
{\bf Abstract} \\[1mm]
\end{center}

   The two--parametric quantum superalgebra $U_{p,q}[gl(2/1)]$ is
consistently defined. A construction procedure for induced representations
of $U_{p,q}[gl(2/1)]$ is described and allows us to 
construct explicitly all (typical and nontypical) finite--dimensional 
representations of this quantum superalgebra. In spite of some
specific features, the present approach is similar to a previously developed method [1]  which, as shown here, is applicable not only to the
one--parametric quantum deformations but also to the
multi--parametric ones.  
\vspace*{1.7cm}   
\begin{center}
\underline{Running title}: Quantum superalgebra $U_{p,q}[gl(2/1)]$
\\[1.7cm] 
MSC--class: 81R50; 17A70

\end{center}
\newpage
\begin{flushleft}
{\Large {\bf I. Introduction}}
\end{flushleft}
\vspace*{2mm}

     In [1], we suggested a method for explicit 
constructions of representations of the one--parametric quantum 
superalgebras  $U_{q}[gl(m/n)]$. When 
applied to $U_{q}[gl[(2/2)]$, this method allowed us to construct 
explicitly all (typical [1] and nontypical [2]) finite--dimensional 
representations of the latter quantum superalgebra. Certainly, as 
emphasized in Refs. 1 and 2, our method is also applicable for other 
quantum superalgebras and we could construct their representations in a similar way. Particularly, we can apply the method to, for example, multiparametric quantum superalgebras, [3--6], etc.. The multiparametric deformations were introduced 
[7] and since considered by a number of authors from different points  of view (see, for example, Refs. 3--15). However, in spite of 
progresses in several aspects (e.g., group--space structures, differential 
calculus, exponential maps, etc.) representation  theory is only 
well developed 
for a few simple cases like  $U_{p,q}[su(2)]$ (see for example Refs. 
8), $U_{p,q}[sl(2/1)]$, [6], etc..
Here, in order to show once again the usefulness of the 
above--mentioned method 
we consider, as a further example, the two-parametric quantum superalgebra 
$U_{p,q}[gl(2/1)]$ 
which, although resembles to the one--parametric quantum superalgebra 
$U_{\sqrt{pq}}[gl(2/1)]$, can not be identified with the latter. 
In this paper we suppose that both $p$ and $q$ are generic, i.e., not 
roots of unity. Following the approach of [1] we can directly 
construct explicit representations of the quantum superalgebra
$U_{p,q}[gl(2/1)]$ induced from some (usually, irreducible)  
finite--dimensional representations of the even subalgebra
$U_{p,q}[gl(2)\oplus gl(1)]$. Since the latter is a stability
subalgebra of $U_{p,q}[gl(2/1)]$ we expect the constructed 
induced representations of $U_{p,q}[gl(2/1)]$ are decomposed into 
finite--dimensional irreducible representations of 
$U_{p,q}[gl(2)\oplus gl(1)]$. For this purpose we shall introduce a 
$U_{p,q}[gl(2/1)]$--basis (i.e., a basis within a 
$U_{p,q}[gl(2/1)]$--module or briefly a basis of $U_{p,q}[gl(2/1)]$) 
which will be convenient for us in investigating the module structure. This basis (see (3.10)) can be expressed in terms of some basis of 
the even subalgebra $U_{p,q}[gl(2)\oplus gl(1)]$ which in turn 
represents a
(tensor) product between a $U_{p,q}[gl(2)]$--basis and a
$gl(1)]$--factor.  
It will be shown that the finite--dimensional representations of
$U_{p,q}[gl(2)]$, i.e., of $U_{p,q}[gl(2)\oplus gl(1)]$ can be
realized in the Gel'fand--Zetlin (GZ) basis.  
The finite--dimensional representations of $U_{p,q}[gl(2/1)]$ 
constructed are irreducible and can be decomposed into
finite--dimensional irreducible representations of the subalgebra
$U_{p,q}[gl(2)\oplus gl(1)]$.\\  

    In section II we shall define the quantum superalgebra  
$U_{p,q}[gl(2/1)]$ and consider how to construct its representations 
induced from representations of the subalgebra $U_{p,q}[gl(2)\oplus 
gl(1)]$. Finite-dimensional representations of $U_{p,q}[gl(2/1)]$ are
constructed in section III where the above--mentioned appropriate basis
is described. The conclusion and some comments are given in section  
IV.\\

   Throughout the paper we shall frequently use the following 
notation 
$$[x]\equiv[x]_{p,q}:={q^{x}-p^{-x}\over q-p^{-1}}\eqno(1.1)$$ for quantum 
deformations of $x$ which are operators or numbers,
$$[X,Y]_{r}:=XY-rYX\eqno(1.2)$$ for r--deformed commutators between two 
operators $X$ and $Y$ and
$$[m]\eqno(1.3)$$ for the highest weights (signatures) of the 
Gel'fand--Zetlin basis vectors $(m)$. We hope this notation will not 
confuse the reader.\\[7mm] 
{\Large {\bf 
II. U$_{p,q}$[gl(2/1)] and its induced representations}}\\

   The two--parametric quantum superalgebra $U_{p,q}[gl(2/1)]$ is
consistently defined through the generators $E_{12}$, 
$E_{21}$, $E_{23}$, $E_{32}$, $E_{ii}$, $i=1,2,3$, and $L$ satisfying 
\begin{tabbing}
\=1234567891234567891\=$[E_{ii},E_{jj}]$1234\= =12\= 01234
\=$[E_{ii},E_{j,j+1}]$\= =
\=$(\delta_{ij}-\delta_{i,j+1})E_{j,j+1}$1\=\kill

~~~~a) the super-commutation relations ($1\leq i,i+1,j,j+1\leq
3$):\\[2mm]
\>\>$[E_{ii},E_{jj}]$\> = \>0,\>\>\>\>(2.1a)\\[1mm]
\>\>$[E_{ii},E_{j,j+1}]$\>=\>$(\delta_{ij}-\delta_{i,j+1})E_{j,j+1}$,\>\>\>\>(2.1b)\\[1mm]
\>\>$[E_{ii},E_{j+1,j}]$\>=\>$(\delta_{i,j+1}-\delta_{ij})E_{j+1,j}$,\>\>\>\>(2.1c)\\[1mm]
\>\>$[L,E_{12}]$\>=\>$[L,E_{21}]$~ = ~$[L,E_{ii}]~=~0$,\>\>\>\>(2.1d)\\[1mm]
\>\>$[E_{12},E_{21}]$\>=\>$\left({q\over p}\right)^{L-h_{1}/2}[h_{1}]$,\>\>\>\>(2.1e)\\[1mm]
\>\>$\{E_{23},E_{32}\}$\>=\>$\left({q\over p}\right)^{-h_{2}}[h_{2}]$.\>\>\>\>(2.1f)\\[1mm]
\>\>$h_{i}$\>=\>$(E_{ii}-{d_{i+1}\over
d_{i}}E_{i+1,i+1}),$\>\>\>\>(2.1g)\\[4mm] with
$d_{1}=d_{2}=-d_{3}=1$,\\[4mm]

~~~~b) the Serre-relations:\\[2mm]
\>\>~~~~~~$E_{23}^{2}$\>=\>~~~~$E_{32}^{2}$\>~~~~~~~=~0,\\[2mm]
\>\>~~$[E_{12},E_{13}]_{p}$\>=\>$[E_{21},E_{31}]_{q}$\>~~~~~~~=~0,
                                        \>\>\>(2.2)\\[2mm]
where\\[2mm]
\>\>~~~~~~~~~~~~$E_{13}$~\>:=\>$[E_{12},E_{23}]_{q^{-1}}$,\\[2mm]
and\\[2mm] 
\>\>~~~~~~~~~~~~$E_{31}$~\>:=\>$-[E_{21},E_{32}]_{p^{-1}}$.\>\>\>\>(2.3) 
\end{tabbing}
are defined as new odd generators which, as we can show, have vanishing squares.
Now the extra--Serre relations are not necessary, unlike in higher rank cases 
[1,2,16]. The commutators between the maximal--spin operator $L$ and the odd 
generators take concrete forms on concrete basis vectors.\\ 

  These generators $E_{ij}$, $i,j= 1,2,3$, are two--parametric deformation 
analogues of the Weyl generators $e_{ij}$ 
$$(e_{ij})_{kl}=\delta_{ik}\delta_{jl},~~i,j,k,l=1,2,3,\eqno(2.4)$$
of the classical (i.e., non--deformed) superalgebra $gl(2/1)$ whose 
universal enveloping algebra $U[gl(2/1)]$ is a classical limit of 
$U_{p,q}[gl(2/1)]$ when $p,q\rightarrow 1$.\\

     From the relations (2.1)--(2.3) we see that every of the odd spaces 
$A_{\pm}$ $$A_{+}= {\normalsize lin.env.}\{E_{13},E_{23}\},\eqno(2.5)$$
$$A_{-}= {\normalsize lin.env.}\{E_{31},E_{32}\},\eqno(2.6)$$
is, as always, a representation space of the even subalgebra 
$U_{p,q}[gl(2/1)_{0}]\equiv U_{p,q}[gl(2)\oplus gl(1)]$ which,  
generated by the generators $E_{12}$, $E_{21}$, $L$ and $E_{ii}$,
$i=1,2,3$, is a stability subalgebra of $U_{p,q}[gl(2/1)]$.
Therefore, we 
can construct representations of $U_{p,q}[gl(2/1)]$ induced from  some 
(finite--dimensional irreducible) representations of $U_{p,q}[gl(2/1)_{0}]$
which are realized in some representation spaces (modules) $V^{p,q}_{0}$ 
being tensor products of $U_{p,q}[gl(2)]$--modules $V^{p,q}_{0, gl_{2}}$
and $gl(1)$--modules (factors) $V^{p,q}_{0,gl_{1}}$. Following [1] 
we demand
$$E_{23}V_{0}^{p,q}=0\eqno(2.7)$$
that is
$$U_{p,q}(A_{+})V_{0}^{p,q}=0.\eqno(2.8)$$
In such a way we turn the $U_{p,q}[gl(2/1)_{0}]$--module $V^{p,q}_{0}$ 
into a $U_{p,q}(B)$--module where
$$B=A_{+}\oplus gl(2)\oplus gl(1).\eqno(2.9)$$
The $U_{p,q}[gl(2/1)]$--module $W^{p,q}$ induced from 
$U_{p,q}[gl(2/1)_{0}]$--module $V^{p,q}_{0}$ is the factor--space
$$W^{p,q}=[U_{p,q}\otimes V_{0}^{p,q}]/I^{p,q}\eqno(2.10)$$
where 
$$U_{p,q}\equiv U_{p,q}[gl(2/1)],\eqno(2.11)$$
while $I^{p,q}$ is the subspace
$$I^{p,q}={\normalsize lin.env.}\{ub\otimes v-u\otimes bv\| u\in U_{p,q},
b\in U_{p,q}(B)\subset U_{p,q}, v\in V_{0}^{p,q}\}.\eqno(2.12)$$
Using the above--given commutation relations (2.1)--(2.2) and the definitions
(2.3) we can prove the following analogue of the Poincar\'e--Birkhoff--Witt 
theorem\\[4mm]
{\it Proposition 1}: The quantum deformation $U_{p,q} := U_{p,q}[gl(2/1)]$ is 
spanned on all possible linear combinations of the elements
$$g = 
(E_{23})^{\eta_{1}}(E_{13})^{\eta_{2}}(E_{31})^{\theta_{1}}
(E_{32})^{\theta_{2}}g_{0},\eqno(2.13)$$
where $\eta_{i}$, $\theta_{i}=0,1$ and $g_{0}\in 
U_{p,q}[gl(2/1)_{0}]\equiv U_{p,q}[gl(2)\oplus gl(1)]$.\\

Then we arrive at the next assertion\\[4mm]
{\it Proposition 2}: The induced $U_{p,q}[gl(2/1)]$--module $W^{p,q}$ is 
the linear span
$$W^{p,q}([m])={\normalsize 
lin.env.}\{(E_{31})^{\theta_{1}}(E_{32})^{\theta_{2}}\otimes v\|v\in 
V_{0}^{p,q},~~\theta_{1}, ~\theta_{2}=0,1\},\eqno(2.14)$$
which is decomposed into (four, at most) finite--dimensional 
irreducible 
modules $V_{k}^{p,q}$ of the even subalgebra $U_{p,q}[gl(2/1)_{0}]$
$$W^{p,q}([m])=\bigoplus _{0\leq k\leq 3}V_{k}^{p,q}([m]_{k}),\eqno(2.15)$$
where $[m]$ and $[m]_{k}$ are some signatures (highest--weights) 
characterizing the module $W^{p,q}\equiv W^{p,q}([m])$ and the modules 
$V_{k}^{p,q}\equiv V_{k}^{p,q}([m]_{k})$, respectively.\\

   As a consequence, for a basis in $W^{p,q}$ we can take all the vectors 
of the form
$$\left |\theta_{1}, \theta_{2}; (m)\right > := 
(E_{31})^{\theta_{1}}(E_{32})^{\theta_{2}}\otimes (m), ~~ \theta_{1},~ 
\theta_{2}=0,1, \eqno(2.16)$$
where $(m)$ is a (GZ, for example,) basis in $V_{0}^{p,q}\equiv 
V_{0}^{p,q}([m])$. We refer to this basis as the induced 
$U_{p,q}[gl(2/1)]$--basis (or simply, the induced basis) in order to 
distinguish 
it from another $U_{p,q}[gl(2/1)]$--basis introduced later and called
a reduced basis.\\

   Any vector $w$ from the module $W^{p,q}$ can be represented as
$$w=u\otimes v,~~~~ u\in U_{p,q},~~~~ v\in V_{0}^{p,q}.\eqno(2.17)$$
Then $W^{p,q}$ is a $U_{p,q}[gl(2/1)]$--module in the sense
$$gw\equiv g(u\otimes v)=gu\otimes v\in W^{p,q}\eqno(2.18)$$
for $g,~u\in U_{p,q}$, $w\in W^{p,q}$ and $v\in V_{0}^{p,q}$.\\[7mm] 
{\Large {\bf
III. Finite--dimensional representations of U$_{p,q}$[gl(2/1)]}}\\

  We can show that finite--dimensional representations of 
$U_{p,q}[gl(2/1)_{0}]$ can be realized in some spaces (modules) $V_{k}^{p,q}$
spanned by the (tensor) products
$$\left[
\begin{array}{lcr}
 
\begin{array}{c}
                        m_{12}~~~m_{22}\\ m_{11}
\end{array}
;
\begin{array}{c}
                        m_{32}=m_{31}\\ m_{31}
\end{array}
\end{array}
\right]
\equiv
\left[
\begin{array}{lcr}
 
\begin{array}{c}
                            [m]_{2}\\ m_{11}
\end{array}
;
\begin{array}{c}
                            [m]_{1}\\ m_{31}
\end{array}
\end{array}
\right]
\equiv
(m)_{gl(2)}\otimes m_{31}\equiv (m)_{k}
\eqno(3.1a)$$
between the (GZ) basis vectors $(m)_{gl(2)}$ of 
$U_{p,q}[gl(2)]$ and the $gl(1)$--factors $m_{31}$, where
$m_{ij}$ are complex numbers such that 
$$m_{12}-m_{11},~ m_{11}-m_{22}\in 
{\bf Z_{+}}\eqno(3.1b)$$ and 
$$m_{32}=m_{31}.\eqno(3.1c)$$
Indeed, any finite--dimensional representation of (not only)
$U_{p,q}[gl(2)]$ is always highest weight and if the generators $L$  
and $E_{ij}$, $i,j=1,2$ are defined on (3.1) as follows
\begin{eqnarray*}
~~~~~~~~~~~~~~~~~ E_{11}(m)_{k}& = &(l_{11}+1)(m)_{k},\\
E_{22}(m)_{k}& = &(l_{12}+l_{22}-l_{11}+2)(m)_{k},\\
E_{12}(m)_{k}& = 
&\left([l_{12}-l_{11}][l_{11}-l_{22}]\right)^{1/2}(m)_{k}^{+11},\\
E_{21}(m)_{k}& =
&\left([l_{12}-l_{11}+1][l_{11}-l_{22}-1]\right)^{1/2}(m)_{k}^{-11},\\
L(m)_{k}&=&{1\over 2}(l_{12}-l_{22}-1)(m)_{k},\\
E_{33}(m)_{k}& = & (l_{31}+1)(m)_{k},
~~~~~~~~~~~~~~~~~~~~~~~~~~~~~~~~~~~~~~~~~~~~~~~~ (3.2a)\\[2mm]
l_{ij}&=&m_{ij}-(i-2\delta_{i,3}),
~~~~~~~~~~~~~~~~~~~~~~~~~~~~~~~~~~~~~~~~~~~~ (3.2b)
\end{eqnarray*}  
where a vector $(m){_k}^{\pm ij}$ is obtained from $(m)$ by replacing  
$m_{ij}$ with $m_{ij}\pm 1$, they really satisfy the commutation 
relations (2.1a)--(2.1e) for $U_{p,q}[gl(2/1)_0]$. The highest weight 
described by the first row (signature) 
$$[m]_{k}=[m_{12},m_{22},m_{32}]\eqno(3.3)$$
of the patterns (3.1) is nothing but an ordered set of the eigen--values 
of the Cartan generators $E_{ii}$, $i=1,2,3$, on the highest weight 
vector $(M)_{k}$ defined as follows
$$E_{12}(M)_{k}=0,\eqno(3.4)$$
$$E_{ii}(M)_{k}=m_{i2}(M)_{k}, \eqno(3.5)$$
The highest weight vector $(M)_{k}$ can be obtained from $(m)_{k}$ by 
setting $m_{11}=m_{12}$
$$(M)_{k}=\left[   
\begin{array}{lcr}
 
\begin{array}{c}
                        m_{12}~~~m_{22}\\ m_{12}
\end{array}
;
\begin{array}{c}
                        m_{32}=m_{31}\\ m_{31}
\end{array}
\end{array}
\right].\eqno(3.6)$$
A lower weight vector $(m)_{k}$ can be derived vice versa from 
$(M)_{k}$ by the formula
\begin{eqnarray*}
~~~~~~~~~(m)_{k}&=&\left({[m_{11}-m_{22}]!\over
[m_{12}-m_{22}]![m_{12}-m_{11}]!}\right)^{1/2}
(E_{21})^{m_{12}-m_{11}}(M)_{k}
~~~~~~~~~~~~~~(3.7)
\end{eqnarray*}\\[2mm]
In particular, for the case $k=0$, instead of the above notations, we
omit the subscript 0, i.e.,  
$$(m)_{0}\equiv (m);~~ [m]_{0}\equiv [m];~~ (M)_{0}\equiv (M),\eqno(3.8)$$
putting
$$m_{i2}=m_{i3},~~~ i=1,2,3,\eqno(3.9)$$
where $m_{i3}$ are some of the complex values of $m_{i2}$,
therefore, $m_{13}-m_{11},~ m_{11}-m_{23}\in {\bf Z_{+}}$.
We emphasize that $[m]$ and $(M)$, because of (2.7), are also,
respectively, the highest weight and the highest weight vector in the
$U_{p,q}[gl(2/1)]$--module $W^{p,q}=W^{p,q}([m])$. Characterizing the
latter module as the whole, $[m]$ and $(M)$ are, respectively,
referred to as the global highest weight and the global highest
weight vector, while $[m]_{k}$ and $(M)_{k}$ are, respectively, the
local highest weights and the local highest weight vectors
characterizing only the submodules $V^{p,q}=V^{p,q}([m]_{k})$.\\

   Following the arguments of [1], for an alternative with (2.16)
basis of $W^{p,q}$ we can choose the union of all the bases (3.1)
which are denoted now by the patterns  
$$\left[
\begin{array}{lcc}
m_{13}& m_{23}& m_{33} \\
m_{12}& m_{22}& m_{32}\\
m_{11}& 0 & m_{31}
\end{array}
\right]_{k}
\equiv
\left[   
\begin{array}{lcr}
 
\begin{array}{c}
                        m_{12}~~~m_{22}\\ m_{11}
\end{array}
;
\begin{array}{c}
                        m_{32}=m_{31}\\ m_{31}
\end{array}
\end{array}
\right]_{k}\equiv (m)_{k}
,\eqno(3.10)$$
where the first row $[m]=[m_{13},m_{23},m_{33}]$ is simultaneously the 
highest weight of the submodule $V^{p,q}=V^{p,q}([m])$ and the whole module 
$W^{p,q}=W^{p,q}([m])$, while the second row $[m]_{k}=[m_{12},m_{22},m_{32}]$
is the local highest weight of some  $U_{p,q}gl[(2/1)_{0}]$--module 
$V^{p,q}_{k}=V^{p,q}_{k}([m]_{k})$ containing the considered vector 
$(m)_{k}$. The basis (3.10) of $W^{p,q}$ is called the 
$U_{p,q}[gl(2/1)]$--reduced basis or simply the reduced basis. The 
latter, as mentioned before and shown later, is convenient for us in 
investigating the module structure of $W^{p,q}$. 

Note once again that the condition 
$$m_{32}=m_{31}\eqno(3.1c)$$ has always to be fulfilled.\\

  The highest weight vectors $(M)_{k}$, now, in the notation (3.10) 
have the form

$$(M)_{k}=\left[
\begin{array}{lcc}
m_{13}& m_{23}& m_{33} \\
m_{12}& m_{22}& m_{32}\\
m_{12}& 0 & m_{31}
\end{array}
\right]_{k}
\eqno(3.11)$$
as for $k=0$ the notation given in (3.8) and (3.9) is also taken into 
account.\\[4mm]
{\it Proposition 3}: The highest weight vectors $(M)_{k}$ are expressed 
in terms of the induced basis (2.16) as follows
\begin{eqnarray*}
~~~~~~~~~~~~~~~ (M)_{0}& = &a_{0}\left|0,0;(M)\right>,~~~~a_{0}\equiv
1,\\[2mm] (M)_{1}& = &a_{1}\left|0,1;(M)\right>,\\[2mm]
(M)_{2}& = &a_{2}\left\{\left|1,0,;(M)\right>
+q^{2l}[2l]^{-1/2}
\left|0,1;(M)^{-11}\right>\right\},\\[2mm]
(M)_{3}& = &a_{3}\left\{\left|1,1;(M)\right>\right\},
~~~~~~~~~~~~~~~~~~~~~~~~~~~~~~~~~~~~~~~~~~~~~~~~~(3.12a)
\end{eqnarray*}
where $a_{i}$, $i=0,1,2,3$, are some numbers depending, in general,
on $p$ and $q$, while $l$ is  
$$l={1\over 2}(m_{13}-m_{23})\eqno(3.12b)$$
Indeed, all the vectors $(M)_{k}$ given above satisfy the condition 
(3.4). From the formulae (3.5) and (3.12) the highest weights
$[m]_{k}$ can be easily identified
\begin{tabbing} \=12345679123456789\= $[m]_{kk}$ \= =x \=
$[m_{13}-1,m_{23}-1,m_{33}+1,m_{43}+1]$,\=\kill
\>\>
$[m]_{0}$ \> = \> $[m_{13}, m_{23}, m_{33}]$,\\[2mm]
\>\>$[m]_{1}$ \> = \> $[m_{13}, m_{23}-1, m_{33}+1]$,\\[2mm] 
\>\>$[m]_{2}$ \> = \> $[m_{13}-1, m_{23},
m_{33}+1]$,\\[2mm] 
\>\>$[m]_{3}$ \> = \> $[m_{13},
m_{23}, m_{33}+2]$~~~~~~~~~~~~~~~~~~~~~~~~~~~~~~~~~~~~~~~(3.13)
\end{tabbing}
\vspace*{2mm}
Using the rule (3.7) we obtain all the basis vectors $(m)_{k}$
\begin{eqnarray*}
~~~~~~~~~~~~~~~ (m)_{0}& \equiv & 
\left[
\begin{array}{lcc}
m_{13}& m_{23}& m_{33} \\
m_{13}& m_{23}& m_{33}\\
m_{11}& 0 & m_{33}
\end{array}
\right]
=\left|0,0,;(m)\right>,\\[4mm]
(m)_{1}&\equiv &
\left[
\begin{array}{lcc}
m_{13}& m_{23}& m_{33} \\
m_{13}& m_{23}-1& m_{33}+1\\
m_{11}& 0 & m_{33}+1
\end{array}
\right] \\[4mm]
& = &a_{1}\left\{-\left(
{[l_{13}-l_{11}]\over
[2l+1]}\right)^{1/2}\left|1,0;(m)^{+11}\right>\right.
\\ &   &\left. +p^{l_{11}-l_{13}}\left(
{[l_{11}-l_{23}]\over
[2l+1]}\right)^{1/2}\left|0,1;(m)\right>\right\}, 
\\[4mm]
(m)_{2}&\equiv &
\left[
\begin{array}{ccc}
m_{13}& m_{23}& m_{33} \\
m_{13}-1& m_{23}& m_{33}+1\\
m_{11}& 0 & m_{33}+1
\end{array}
\right] \\[4mm]
& = &a_{2}\left\{\left({q\over p}\right)^{l_{13}-l_{11}-1}\left(
{[l_{11}-l_{23}]\over
[2l]}\right)^{1/2}\left|1,0;(m)^{+11}\right>\right.
\\ &   &\left. +q^{l_{13}-l_{23}-1}p^{l_{11}-l_{13}+1}\left(
{[l_{13}-l_{11}]\over
[2l]}\right)^{1/2}\left|0,1;(m)\right>\right\},
\\[4mm]
(m)_{3}&\equiv &
\left[
\begin{array}{ccc}
m_{13}& m_{23}& m_{33} \\
m_{13}-1& m_{23}-1& m_{33}+2\\
m_{11}& 0 & m_{33}+2
\end{array}
\right] \\[4mm]
&=& a_{3}\left|1,1;(m)\right>,
~~~~~~~~~~~~~~~~~~~~~~~~~~~~~~~~~~~~~~~~~~~~~~~~~~~~~~(3.14)
\end{eqnarray*}
where  $l_{ij}$ and $l$ are given in (3.2b) and (3.12b),
respectively. Here,  we omit the subscript $k$ in the above patterns
since there is no degeneration between them. 
The formulae (3.14), in fact, represent the way in which the reduced
basis (3.10) is written in terms of the induced basis (2.16). From
(3.14) we can derive their inverse relation 
\begin{eqnarray*}
~~~~~~~~~ \left|1,0;(m)\right>& = &(m)\\[2mm]
\left|1,0;(m)\right>& = &-{1\over a_{1}}q^{l_{11}-l_{23}-1}
\left({[l_{13}-l_{11}+1]\over [2l+1]}
\right)^{1/2}(m)_{1}^{-11}\\[2mm] &   & +{1\over
a_{2}p}q^{l_{11}-l_{13}} {\left([l_{11}-l_{23}-1][2l] \right)^{1/2}\over
[2l+1]}(m)_{2}^{-11},\\[2mm] 
\left|0,1;(m)\right>& = &{1\over a_{1}}
\left({[l_{11}-l_{23}]\over [2l+1]}   
\right)^{1/2}(m)_{1}\\[2mm] &   & + {1\over a_{2}}\left({p\over 
q}\right)^{l_{13}-l_{11}-1} {\left([l_{13}-l_{11}][2l]\right)^{1/2}\over 
[2l+1]} (m)_{2},\\[2mm] 
\left|1,1;(m)\right>& = &{1\over 
c_{3}}(m)^{-11}_{3}.
~~~~~~~~~~~~~~~~~~~~~~~~~~~~~~~~~~~~~~~~~~~~~~~~~~~~~~~~~~~~(3.15)
\end{eqnarray*}
\vspace*{2mm}
    Now we are ready to compute all the matrix elements of the
generators in the basis (3.10). As we shall see, the latter
basis allows a clear descriptios of a decomposition of a
$U_{p,q}[gl(2/1)]$--module 
$W^{p,q}$ in irreducible $U_{p,q}[gl(2/1)_{0}]$--modules
$V^{p,q}_{k}$. Since the finite--dimensional 
representations of the $U_{p,q}[gl(2/1)]$ in some basis are
completely defined by the actions of the even generators and the 
odd Weyl--Chevalley ones $E_{23}$ and $E_{32}$ in the same basis, 
it is sufficient to write down the matrix elements of these generators
only. For the even generators the matrix elements have already been
given in (3.2), while for $E_{23}$ and $E_{32}$, using the relations
(2.1)--(2.3), (3.14) and (3.15) we have 
\begin{eqnarray*}
E_{23}(m) &=& 0,\\[2mm]
E_{23}(m)_{1} &=&a_{1}\left({p\over q}\right)^{l_{23}+l_{33}+3}
\left({[l_{11}-l_{23}]\over [2l+1]}\right)^{1/2}[l_{23}+l_{33}+3] 
(m),\\[4mm]
E_{23}(m)_{2} &=&a_{2}\left({p\over q}\right)^{l_{23}+l_{33}+4}
\left({[l_{13}-l_{11}]\over [2l]}\right)^{1/2}[l_{13}+l_{33}+3]
(m),\\[4mm]
E_{23}(m)_{3} &=&a_{3}\left({p\over 
q}\right)^{l_{13}+l_{23}+l_{33}-l_{11}+2}
\left\{{1\over a_{1}q}\left({[l_{13}-l_{11}]\over 
[2l+1]}\right)^{1/2}[l_{13}+l_{33}+3](m)_{1}\right. ,\\[2mm]
&&\left. 
-{1\over a_{2}p}\left({[l_{11}-l_{23}] 
[2l]}\right)^{1/2}{[l_{23}+l_{33}+3]\over [2l+1]}(m)_{2}\right\}
~~~~~~~~~~~~~~~~~~~~~~~~~(3.16a)
\end{eqnarray*}
and
\begin{eqnarray*}
~~~~~~~~~~~ E_{32}(m) &=& {1\over a_{1}}\left({[l_{11}-l_{23}]\over 
[2l+1]}\right)^{1/2}(m)_{1},\\[2mm]
&&+{1\over a_{2}}\left({p\over q}\right)^{l_{13}-l_{11}-1}
{\left([l_{13}-l_{11}][2l]\right)^{1/2}\over 
[2l+1]}(m)_{2}\\[4mm] 
E_{32}(m)_{1} &=&{a_{1}\over 
a_{3}}p
\left({[l_{13}-l_{11}]\over [2l+1]}\right)^{1/2}
(m)_{3},\\[4mm]
E_{32}(m)_{(2)} &=&-{a_{2}\over a_{3}}p\left({q\over 
p}\right)^{l_{13}-l_{11}-1}
\left({[l_{11}-l_{23}]\over [2l]}\right)^{1/2}
(m)_{3},\\[4mm]
E_{32}(m)_{3} &=&0.
~~~~~~~~~~~~~~~~~~~~~~~~~~~~~~~~~~~~~~~~~~~~~~~~~~~~~~~~~~~~~~~~~~~~ (3.16b)
\end{eqnarray*}
{\it Proposition 4}: The finite--dimensional representations (3.16) of 
$U_{p,q}[gl(2/1)]$ are irreducible and called typical if only if the 
condition
$$[l_{13}+l_{33}+3][l_{23}+l_{33}+3]\neq 0\eqno(3.17)$$
holds.\\

   When this condition (3.17) is violated, i.e. one  of the following
condition pairs 
$$[l_{13}+l_{33}+3]=0\eqno(3.18a)$$
and 
$$[l_{23}+l_{33}+3]\neq 0\eqno(3.18b)$$
or
$$[l_{13}+l_{33}+3]\neq 0\eqno(3.19a)$$
and
$$[l_{23}+l_{33}+3]=0\eqno(3.19b)$$
(but not both (3.18a) and (3.19b) simultaneously) holds,
the module $W^{p,q}$ is no longer irreducible but indecomposable. 
However, there exists an invariant 
subspace, say $I_{k}^{p,q}$, of $W^{p,q}$ such that the 
factor--representation in the factor--module
$$W_{k}^{p,q}:=W^{p,q}/I_{k}^{p,q}\eqno(3.20)$$
is irreducible. We say that is a nontypical representation in a 
nontypical module $W_{k}^{p,q}$. Then, as in [2], it is not difficult 
for us to prove the following assertions\\[4mm]
{\it Proposition 5}: 
$$V_{3}^{p,q}\subset I_{k}^{p,q},\eqno(3.21)$$
and
$$V_{0}^{p,q}\cap I_{k}^{p,q}=\O.\eqno(3.22)$$
From (3.16)--(3.18) we can 
easily find all nontypical representations of 
$U_{p,q}[gl(2/1)]$ which are classified into two classes.\\[4mm]

 {\bf III.1. Class 1 nontypical representations:}\\

  This class is characterized by the conditions $(3.18a)$ and 
$(3.18b)$ which for generic $p$ and $q$ take the forms 
$$l_{13}+l_{33}+3=0,\eqno(3.18a')$$
and 
$$l_{23}+l_{33}+3\neq 0,\eqno(3.18b')$$
respectively.
In other words, we have to replace everywhere all $m_{33}$ by $-m_{13}-1$
and keep $(3.18b')$ valid. Thus we have\\[4mm]
{\it Proposition 6}:
$$I_{1}^{p,q}=V_{3}^{p,q}\oplus V_{2}^{p,q}.\eqno(3.23)$$

   Then the class 1 nontypical representations in
$$W_{1}^{p,q}=W_{1}^{p,q}([m_{13},m_{23},-m_{13}-1])\eqno(3.24)$$
are given through (3.16) by keeping the conditions (3.18) (i.e.,
$(3.18a')$ and $(3.18b')$) and  
replacing all vectors belonging to $I_{1}^{p,q}$ with 0:
\begin{eqnarray*}
~~~~~~~~~ E_{23}(m)&=&0,\\[2mm]
E_{23}(m)_{1}&=&a_{1}\left({p\over q}\right)^{l_{23}-l_{13}}
\left({[l_{11}-l_{23}]\over [2l+1]}\right)^{1/2}[l_{23}-l_{13}](m)
~~~~~~~~~~~~~~~~~~~~(3.25a)
\end{eqnarray*}
and
\begin{eqnarray*}
~~~~~~~~~~~~~~~~~~~~~~~~ E_{32}(m)&=&{1\over a_{1}}\left({[l_{11}-l_{23}]\over 
[2l+1]}\right)^{1/2}(m)_{1}\\[2mm]
E_{32}(m)_{1}&=&0.~~~~~~~~~~~~~~~~~~~~~~
~~~~~~~~~~~~~~~~~~~~~~~~~~~~~~~~~~~(3.25b)
\end{eqnarray*}
\vspace*{2mm}

  {\bf III.2. Class 2 nontypical representations:}\\

   For this class nontypical representations we must keep the conditions

$$l_{13}+l_{33}+3\neq 0,\eqno(3.19a')$$
and
$$l_{23}+l_{33}+3= 0.\eqno(3.19b')$$
derived respectively from $(3.19a)$ and $(3.19b)$ when the 
deformation parameters $p$ and $q$ are generic. 
Equivalently, we have to replace everywhere all $m_{33}$ by $-m_{23}$
and keep $(3.19a')$ valid.\\

   Now the invariant subspace $I_{2}^{p,q}$ is the following\\[4mm]
{\it Propositions 7}:
$$I_{2}^{p,q}=V_{3}^{p,q}\oplus V_{1}^{p,q}.\eqno(3.26)$$

  The class 2 nontypical representations in 

$$W_{2}^{p,q}=W_{2}^{p,q}([m_{13},m_{23},-m_{23}])\eqno(3.27)$$
are also given through (3.16) but by keeping the conditions (3.19)
(i.e., $(3.19a')$ and $(3.19b')$)
valid and replacing all vectors belonging to the invariant by 0
subspace $I_{2}^{p,q}$:

\begin{eqnarray*}
~~~~~~~~~~~~~~ E_{23}(m)&=&0,\\[2mm]
E_{23}(m)_{2}&=&a_{1}{p\over q}
\left({[l_{13}-l_{11}]\over [2l]}\right)^{1/2}[2l+1](m)
~~~~~~~~~~~~~~~~~~~~~~~~~~~~~(3.28a)
\end{eqnarray*}
and
\begin{eqnarray*}
~~~~~~~~~~~~~~~~~~~~~~~~~~ E_{32}(m)&=&{1\over 
a_{2}}{\left([l_{13}-l_{11}][2l]\right)^{1/2}\over
[2l+1]}(m)_{2}\\[2mm]
E_{32}(m)_{2}&=&0. ~~~~~~~~~~~~~~~~~~~~~~~~~~~~~~~~~~~~~~~~~~~~~~~~~~~~~~~ (3.28b)
\end{eqnarray*}

  In order to complete this section we emphasize that nontypical
representations have only been well investigated for a  
few cases of both classical and quantum superalgebras (see, in this
context, the Conclusion in Ref. 2 and also some comments in Ref. 17).
Therefore, the present results can be considered as a small
step forward to this direction.\\[1.3cm]
{\Large {\bf IV. Conclusion}}\\

   We have just defined the two--parametric quantum superalgebra 
$U_{p,q}[gl(2/1)]$ and constructed at generic deformation parameters
all its  
typical and nontypical representations leaving the coefficients
$a_{i}$, $i=1,2,3$, as free parameters which can be fixed by some
additional conditions, for example, the hermiticity condition.
As an intermediate step (which, however, is of independent
interest) we also introduced the reduced basis (3.10) which, 
as it is an extension of the Gel'fand--Zetlin basis to the present 
case, 
is appropriate for a clear description of decompostions of
$U_{p,q}[gl(2/1)]$--modules into irreducible 
$U_{p,q}[gl(2/1)_{0}]$--modules. Although the present approach has 
some specific features it is 
similar to the one in Ref. 1. That shows once again the usefulness 
of the method of Ref. 1 which is thus applicable not only 
to the one--parametric quantum deformations but also to the
multi--parametric ones.\\ 

 As the general procedure has been given, the next step is to consider the case of non--generic $p$ and $q$ or to construct 
representations of larger quantum superalgebras like $U_{p,q}[gl(n/1)]$, 
$U_{p,q}[gl(n/m)]$, etc. for both generic and non--generic deformation 
parameters. Let us emphasize once again that our approach avoids the
use of the Clebsch--Gordan coefficients which are not always known,
especially for higher rank (classical and quantum) algebras and
multi--parametric deformations.\\[9mm] 
{\large {\bf Acknowledgements}}\\

   Recently, I learnt that Professor Asim Orhan Barut is no longer 
amongst us. I should like to dedicate this paper to Professor Barut -- an 
outstanding scientist and man. His book (with R. Raczka) [18] on the 
group representation theory and its applications is one of the best
books I ever had and is frequently used.\\
   
   I would like to thank Professor G. Altarelli and Professor G. Veneziano
for the kind hospitality at the Theory Division, CERN, Geneva, Switzerland.
I am thankful to Professor T. Palev for bringing Ref. 6 to my
attention after he had read the present paper.  
\newpage
\begin{flushleft}
{\Large {\bf References}}
\end{flushleft}
\begin{enumerate}
\item Nguyen Anh Ky, J. Math. Phys. {\bf 35}, 2583 (1994) or hep--th/9305183.
\item Nguyen Anh Ky and N.I. Stoilova, J. Math. Phys. {\bf 36}, 5979
(1995) and hep--th/9411098.
\item Yu. Manin, Comm. Math. Phys., {\bf 123}, 169 (1989).
\item R. Chakrabarti and R. Jagannathan, {\it On Hopf structure of
$U_{p,q}(gl(1|1))$ and the universal ${\cal T}$--matrix of
$Fun_{p,q}(Gl(1|1))$}, ICTP preprint IC/94/254 or  
hep--th/9409161; Z. Phys. C: Part. \& Fields, {\bf 66}, 523 (1991).
\item L. Dabrowski and Lu-uy Wang, Phys. Lett. {\bf B 266}, 51 (1991).
\item R. Zhang, J. Phys. A: Math. Gen., {\bf 23}, 817 (1994).
\item Yu. Manin, {\it Quantum groups and non--commutative geometry}, 
Centre des Recherchers Math\'ematiques, Montr\'eal (1988).
\item M. Kibler, Mod. Phys. Lett. A {\bf 10}, 51 (1995); J. Phys. G: {\bf
20}, L13 (1994) and references therein.
\item N. Reshetikhin, Lett. math. Phys. {\bf 20}, 331 (1990).
\item A. Sudbery, J.Phys. A: Math. Gen., {\bf 23}, L697 (1990).
\item A. Schirmacher, J. Wess and B. Zumino, Z. Phys. C: Part. \& Fields, 
{\bf 49}, 317 (1991).
\item A. Schirmacher, Z. Phys. C: Part. \& Fields, {\bf 50}, 321 (1991).
\item R. Jagannathan and J. Van der Jeugt, J.Phys. A: Math. Gen., {\bf 
28}, 2819 (1995); {\it The exponential map for representations of $U_{p,q}(gl(2))$}, q--alg/9507009.
\item C. Fronsdal and A. Galino, Lett. Math. Phys., {\bf 27}, 59 (1993).
\item C. Fronsdal and A. Galino, Lett. Math. Phys., {\bf 34}, 25 (1995).
\item R. Floreanini, D. Leites and L. Vinet, Lett. Math. Phys.
{\bf 23}, 127 (1991); 
M. Scheunert, Lett. Math. Phys. {\bf 24}, 173 (1992); S. M. Khoroshkin and
V. N.
Tolstoy, Comm. Math. Phys. {\bf 141}, 599 (1991).
\item B. Abdesselam, D. Arnaudon and A. Chakrabarti, {\it Representations of $U_{q}[sl(n)]$ at roots of unity}, preprint 
ENSLAPP--A--506/95 or q-alg/9504006.
\item A. O. Barut and R. Raczka, {\it Theory of Group Representations and Applications} (Polish Scientific, Warsaw, 1977); this book was later published in several other English and Russian editions.
\end{enumerate}   
\end{document}